\documentclass[11pt,twoside,a4paper]{amsart}
\usepackage[french]{babel}
\usepackage{amssymb,bbm}
\usepackage{amsmath}
\usepackage{amsthm}
\usepackage{hyperref}
\usepackage[dvips]{graphics}
\usepackage{graphicx}
\usepackage{latexsym}
\usepackage[normal]{subfigure}
\input{epsf.sty}
\include{psfig}
\usepackage[all]{xy}
\usepackage{eufrak}
\usepackage{mathrsfs}
\usepackage{lipsum}
\usepackage{mdframed}
\usepackage{perpage}
\usepackage{amsfonts}
\usepackage{url,color,srcltx}
\begin{document}
\def\K{\mathbb{K}}
\def\R{\mathbb{R}}
\def\C{\mathbb{C}}
\def\Z{\mathbb{Z}}
\def\Q{\mathbb{Q}}
\def\D{\mathbb{D}}
\def\N{\mathbb{N}}
\def\T{\mathbb{T}}
\def\P{\mathbb{P}}
\def\A{\mathscr{A}}
\def\CC{\mathscr{C}}
\renewcommand{\theequation}{\thesection.\arabic{equation}}
\newtheorem{theorem}{Theorem}[section]
\newtheorem{lemma}{Lemma}[section]
\newtheorem{corollary}{Corollaire}[section]
\newtheorem{prop}{Proposition}[section]
\newtheorem{definition}{Definition}[section]
\newtheorem*{preuve}{Proof}
\newtheorem{remark}{RemarK}[section]
\newtheorem{Example}{Example}[section]
\newtheorem{notation}{Notation}
\newtheorem{Non}{Cons\'equence}
\bibliographystyle{plain}
                                    
\title[The $\partial \bar{\partial}$-problem for a forms with boundary value in currents sense ~~]{\textbf {The $\partial \bar{\partial}$-problem for a differential forms with boundary value in currents sense defined in a contractible completely strictly pseudoconvex domain of a complex manifold } }
\author[ S.\  Sambou   \& S.\  Sambou ]
{Salomon Sambou \&  Souhaibou Sambou}
\address{
Mathematics Department\\UFR of Science and Technology \\ Assane Seck University of Ziguinchor, BP: 523 (S{\'e}n{\'e}gal)}

\email{ssambou@univ-zig.sn  \& s.sambou1440@zig.univ.sn  }

   


\maketitle
\begin{abstract}
We solve the $\partial \bar{\partial}$-problem for the differential forms of class $C^\infty$ with boundary value in currents sense defined on a contractible completely strictly pseudoconvex domain of  a complex manifold.
\end{abstract}
\vskip 2mm
\noindent
\maketitle  
\section*{Introduction}
Let $M$ be a differentiable manifold and $\Omega \subset M$ be a domain. In this work,
we try to solve the $\partial \bar{\partial}$-equation  for the differential forms of class $C^\infty$ with boundary value in currents sense.\\
We follow the same steps of resolution of $[7]$ where the same problem has been solved for open sets of $\mathbb{C}^n$. The convex part of $[7]$ can be seen as the local analogues of the convex part of this work. For the concave part, inspired by lemma $4.3$ of $[9]$ and the concave case of $[7]$, we obtain the
local results also. It is known that we do go through recollement from local results to global results, the partition of the unit not being holomorphic. We will move from local to global by using the classic sheaves theory. The result obtained in this direction is the following :
\begin{theorem} \label{zhou}
Let $M$ be a complex analytic manifold of dimension $n$ and $\Omega \subset \subset M$ a contractible domain that is completely strictly pseudoconvex with smooth boundary of class $C^\infty$. We Suppose that $H^j(b\Omega)$ is trivial for $1 \leq j \leq 2n-2$. Then the equation $\partial \bar{\partial}u = f$, where $f$ is a $(p, q)$-form of class $C^\infty$, $d$-closed with boundary value in currents sense for $1 \leq p \leq n-1$ and $1 \leq q \leq n-1$ admits a solution $u$ which is a $(p-1, q-1)$-form of class $C^\infty$ with boundary value in currents sense.
\end{theorem}
We also consider the concave version of the previous theorem and
we obtain :
\begin{theorem} \label{aff}
 Let $M$ be a complex analytic manifold of dimension $n$ and $D \subset \subset M$ a contractible domain that is completely strictly pseudoconvex with smooth boundary of class $C^\infty$. We Suppose that $M$ is a $(n-1)$-convex extension of $D$ and it is also a contractible extension of $D$ with $H^j(b\Omega)$ is trivial for $1 \leq j \leq 2n-2$. Let $\Omega = M \setminus \bar{D}$. If $\stackrel{\circ}{\bar{\Omega}} = \Omega$ and $f$ is a $(p, q)$-form of class $C^\infty$, $d$-closed with boundary value in currents sense on $\Omega$ for $1 \leq p \leq n-1$ and $1 \leq q \leq n-1$, then there exists
a $(p - 1, q - 1)$-form $u$ of class $C^\infty$ with boundary value in currents sense on $\Omega$ such that $\partial \bar{\partial}u = f$.
\end{theorem}
\section{Preliminaries}
\begin{definition}{(voir $[2]$)}
Let $X$ be a differentiable manifold and $\Omega \subset \subset X$ be a contractible domain. We say that $X$ is an contractible extension of $\Omega$, if there is an exhaustive  sequence $(\Omega_n)_{n \in \mathbb{N}}$ of contractibles domains  such as\\
\begin{center}
$\forall~~ n \in \mathbb{N}$,  $\Omega \subset \subset \Omega_n \subset \subset X$.
\end{center}
\end{definition}
\begin{Example}
When $X = \mathbb{C}^n$, then $\mathbb{C}^n$ is a contractible extension of the unit ball $B$.
\end{Example}
\begin{definition}
Let $X$ be a complex analytic  manifold of dimension $n$.
\begin{enumerate}
\item[(1)] A function $\rho$ of class $C^\infty$ on $X$ is called $n$-convex (respectively
$n$-concave) if its form of Levi admits $n$ eigenvalues ​​strictly
positives (respectively strictly negatives).
\item[(2)] Let $\Omega \subset \subset X$ be a relatively compact domain of $X$. $\Omega$ is said to be completely strictly pseudoconvex if there is a function $n$-convex $\varphi$ defined in a neighborhood $U_{\bar{\Omega}}$ of $\Omega$ such as
 $\Omega = \{ z \in U_{\bar{\Omega}}~\vert~ \varphi(z) < 0 \}$.
\item[(3)] $X$ is an $(n - 1)$-concave extension of $\Omega$ if:
\begin{enumerate}
\item[(i)] $\Omega$ meets all the connected components of $X$.
\item[(ii)] There exists an $n$-concave function $\varphi$ defined on a neighborhood $U$ of
$X \setminus \Omega$ such that $\Omega \cap U = \{ z \in U~\vert~ \varphi(z) <0\}$ and for any real $\alpha$ with $0< \alpha < \sup_{z \in U} \varphi(z)$ the set $\{z \in U~\vert~0 \leq \varphi(z) \leq \alpha\}$ is compact.
\end{enumerate}
\end{enumerate}
\end{definition}
\begin{definition}
Let $M$ be a differentiable manifold and $\Omega \subset \subset X$ be a smooth domain with boundary of class $C^\infty$ of defining function $\rho$. Let\\
$\Omega_\varepsilon = \{z \in \Omega ~\vert~ \rho(z)<-\varepsilon\}$ and $b\Omega_\varepsilon$ the boundary of $\Omega_\varepsilon$.\\
Let $f$ be a function of class $C^\infty$ on $\Omega$. It is said that $f$ has a boundary value in distributions sense, if there is a distribution $T$ defined on the
boundary $b\Omega$ of $\Omega$ such that for any function $\varphi \in C^\infty (b\Omega)$, we have:\\
\[ \lim_{\varepsilon \rightarrow 0} \int_{b\Omega_\varepsilon} f \varphi_\varepsilon d\sigma = <T,\varphi>\]
where $\varphi_\varepsilon = i_{\varepsilon}^* \tilde{\varphi}$  with $\tilde{\varphi}$ an extension of $\varphi$ to $\Omega$ and $i_{\varepsilon} : b\Omega_\varepsilon  \rightarrow b\Omega$  the injection
canonical; $d\sigma$ denotes the element of volume.\\
A differential form of class $C^\infty$ on $\Omega$ admits a boundary value in currents sense if its coefficients have a boundary value in distributions sense.
\end{definition}
\begin{definition}
We say that a function $f$ of class $C^\infty$ defined on $\Omega$ is of polynomial growth of order $N \geq 0$, if there is a constant $C$ such that for all $z \in \Omega$, we have:
\begin{center}
$\vert f(z) \vert \leq \frac{C}{d(z)^N}$
\end{center}
where $d(z)$ is the distance from $z$ to the boundary of $\Omega$.
\end{definition}
\begin{definition}
Let $\Omega \subset M$ be a domain in a differentiable manifold $M$. The current T defined on $\Omega$ is said to be extensible if $T$ is the restriction to $\Omega$ of a current $\tilde{T}$ defined on $M$.\\
It is known from [3] that if $\stackrel{\circ}{\bar{\Omega}} = \Omega$, then the extensible currents are
the topological duals of the differential forms with compact support in $\bar{\Omega}$.\\
We will have in all the following to consider the domains $\Omega$ verifying $\stackrel{\circ}{\bar{\Omega}} = \Omega$.
\end{definition}
\begin{notation}
Let $\Omega$ be a domain in a complex analytic manifold $M$ of dimension $n$.\\
We denote by $O_\Omega$ the sheaf of holomorphic functions on $\Omega$, $\check{O}_\Omega$ the one on $\bar{\Omega}$ of germs of holomorphic functions on $\Omega$ with boundary value in distributions sense and $\mathcal{F}^{0,r}(\Omega)$ the one on $\bar{\Omega}$ of the $(0, r)$-forms with boundary value in currents sense.\\
We notice
$\check{H}^r(\Omega)$ the $r^{th}$ De Rham cohomology group of extensible currents defined on $\Omega$, $H^r(\Omega)$ the $r^{th}$ De Rham cohomology group of
differentiable forms of class $C^\infty$ defined on $\Omega$ and $H^r(b\Omega)$ the
$r^{th}$ De Rham cohomology group of differentiable forms of class $C^\infty$ defined on $b\Omega$. The $r^{th}$ group of De Rham cohomology of differentiable forms of class $C^\infty$ with boundary value in currents sense on $\Omega$ is noted $\tilde{H}^r(\Omega)$.\\
We denote by $H^r(\Omega,O_\Omega)$ (respectively $H^r(\Omega,\check{O}_\Omega)$) the $r^{th}$ $\check{C}$ech cohomology group of differentiable forms  defined on $\Omega$ with value in the
sheaf $O_\Omega$ (respectively $\check{O}_\Omega$).
\end{notation}
\section{Solving the equation $du=f$}
\begin{theorem} \label{1}
Let $M$ be a differentiable manifold of dimension $n$ and $\Omega \subset \subset M$ be a relatively compact domain with smooth boundary of class $C^\infty$. We suppose that $\Omega$ is contractible and $H^j(b\Omega)$ is trivial for $1 \leq j \leq n - 2$. Then the equation $du=f$ where $f$ is a $r$-form of class $C^\infty$ defined on $\Omega$ with boundary value in currents sense and $d$-closed admits a solution $u$ which is a $(r - 1)$-form of class $C^\infty$ defined on $\Omega$ with boundary value in currents sense for $1 \leq r \leq n - 1$.
\end{theorem}
\begin{preuve} 
According to $[5]$ if $f$ is a $d$-closed differential form with boundary value in currents sense on $\Omega$ then $[f]$ is a extensible current. Since $\bar{\Omega}$ is compact, the current $[f]$ is of finite order. According to $[1]$, $\check{H}^r(\Omega)= 0$, there exists an extensible current $S$ defined on $\Omega$ such that $dS = f$. Let $\tilde{S}$ be an extension with compact support in $\bar{\Omega}$ of $S$. according to ($[4]$ page $40$)
\begin{center}
$\tilde{S} = R\tilde{S} + Ad\tilde{S} + dA\tilde{S}$
\end{center}
$d\tilde{S} =d( R\tilde{S} + Ad\tilde{S})$ and $( R\tilde{S} + Ad\tilde{S})_{\vert \Omega}$ is another solution of $du = f$ on $\Omega$ and it is an extensible current. $R\tilde{S}$ is an differential form of class $C^\infty$ on $\Omega$ and $A$ does not increase the singular support of $d\tilde{S}$. Since $d \tilde{S}_{\vert \Omega}= f$ which is of class $C^\infty$ then $Ad\tilde{S}_{\vert \Omega}$ is of class $C^\infty$. Thus $( R\tilde{S} + Ad\tilde{S})_{\vert \Omega}$ is of class $C^\infty$. It remains to show that $( R\tilde{S} + Ad\tilde{S})_{\vert \Omega}$ admits a boundary value in currents sense. $R\tilde{S}$ is an differential form of class $C^\infty$ therefore admits a boundary value in currents sense.\\
Let  $(\varphi_j)_{j \in J}$ be a partition of the unit subordinated to a finite recovery $(U_j)_{j \in J}$ of $\bar{\Omega}$ by the open set of local coordinates.\\
We have $Ad \tilde{S} = \displaystyle {\sum_{j \in J} A \varphi_jd \tilde{S} }$ with $A\varphi_j d \tilde{S}$ have compact support in $U_ j$.
\begin{center}
$Ad \tilde{S}_{\vert \Omega} = \displaystyle {\sum_{j \in J} (A \varphi_j d \tilde{S})_{\vert \Omega}}$.
\end{center}
If $U_j \subset \Omega$, then $ A \varphi_j d \tilde{S}$ is of class $C^\infty$ with compact support in $\Omega$, so $(A \varphi_j d \tilde{S})_{\vert \Omega}$ admits a boundary value in currents sense.\\
If $U_j \nsubseteq \Omega$ and $U_j \cap b \Omega \neq \emptyset$; let us show that $ A \varphi_j d \tilde{S}$ admits a boundary value in currents sense.\\
Since $\varphi_j$ is with support in $U_j$ which is an open set of coordinate, so
we are reduced to a bounded domain of $\mathbb{R}^n$. $ A \varphi_j d \tilde{S}$  is of the same nature
than the action of the newton kernel $E(x, y)$ on $\varphi_j d \tilde{S}$. Since $d \tilde{S}$ is order
$l$ with compact support in $\bar{\Omega}$ which prolongs $f$, thus according to $[7]$, $ A \varphi_j d \tilde{S}$ admits a boundary value in currents sense.\\
Since $Card(J)< \infty$, so $\displaystyle {\sum_{j \in J} (A \varphi_j d \tilde{S})_{\vert \Omega}}$ admits a boundary value in currents sense. Thus we have
\begin{center}
$\tilde{H}^{r}(\Omega) = 0$.
\end{center}
\end{preuve}
\section{ Solving $\partial \bar{\partial}$ for a differential forms with boundary value in currents sense}
As consequences of Theorem $4.1$ of $[7]$, we have the following corollary:
\begin{corollary} \label{4}
Let $X$ be a complex analytic manifold of dimension $n$ and $\Omega$ be a completely strictly $(n-1)$-convex  domain of $X$. Let $f$ be a differential form of bidegre $(0, r)$ of class $C^\infty$ with boundary value in currents sense, defined on $\Omega$ and $\bar{\partial}$-closed for $1 \leq r \leq n - 1$. Then there exists a $(0, r - 1)$-form $g$ of class $C^\infty$ defined on $\Omega$ with boundary value in currents sense such that $\bar{\partial}g = f$.
\end{corollary}
\begin{preuve}
Consider the following sequence:
\begin{equation} 
0 \rightarrow \check{O}_\Omega \rightarrow \mathcal{F}^{0,0}(\Omega) \rightarrow \mathcal{F}^{0,1}(\Omega) \rightarrow  \cdots \rightarrow \mathcal{F}^{0,n-1}(\Omega) \rightarrow \mathcal{F}^{0,n}(\Omega) \rightarrow 0.
\end{equation}
So we have a complex $(\mathcal{F}^{0,\bullet}(\Omega), \bar{\partial})$ of $\bar{\partial}$ for differential forms defined on $\Omega$ with boundary value in currents sense. Thanks to the local solving of $\bar{\partial}$  in $[7]$, this complex is an acyclic solving of sheaf $\check{O}_\Omega$. This leads to the isomorphism $H^r(\Omega,\check{O}_\Omega)\simeq \tilde{H}^{0,r}(\Omega)$. According to Corollary $5.2$ of $[8]$ we have $H^r(\Omega,\check{O}_\Omega)=0$.\\ 
So
\begin{center}
$\tilde{H}^{0,r}(\Omega) = 0$.
\end{center}
\end{preuve}
More generally, if $0 \leq q \leq n - 1$, we have the following theorem:
\begin{theorem} \label{q}
Let $X$ be a complex analytic manifold of dimension $n$ and $\Omega \subset \subset X$ be a domain with smooth boundary of class $C^\infty$ and completely strictly $q$-convex for $0 \leq q \leq n - 1$. Let $f$ be a $(n, r)$-form of class $C^\infty$, defined on $\Omega$ and $\bar{\partial}$-closed with boundary value in currents sense for $1 \leq n - q \leq r \leq n$. Then there exists a $(n, r - 1)$-form $g$ of class $C^\infty$ defined on $\Omega$ with boundary value in currents sense such as $\bar{\partial}g = f$.
\end{theorem}
\begin{preuve}
Since we consider the $(n, r)$-forms then the operator $\bar{\partial}$ is equal to the exterior differentiation operator d.\\
Let $f \in \tilde{H}^{n,r}(\Omega)$, $[f]$ is a extensible current defined on $\Omega$. Since
According to $[8]$ we have $\check{H}^{n,r}(\Omega) = 0$ for $n-q \leq r \leq n$ then there exists a $(n, r-1)$-extensible current $u$ defined on $\Omega$ such that $du = f$. Let $S$ be an extension of $u$ with compact support in $\bar{\Omega}$. $S$ is of finite order, $F = dS$ is
a extensible current of finite order noted $m$ and $F_{\vert \Omega} = f$. According to ($[4]$ page
$40$), we have
\begin{center}
$S = RS + AdS + dAS$.
\end{center}
Now $dS= F$ $\Rightarrow$ $S= RS+ AF+ dAS$ then $dS= d(RS+ AF)= F$ so $(RS + AF)_{\vert \Omega}$ is another solution of the equation of $du = f$. Now $RS$ is a Differential form of class $C^\infty$ with compact support therefore admits a boundary value in currents sense. Since $A$ does not increase the singular support, we have if $F$ is of class $C^\infty$ on $\Omega$ then $AF$ is also of class $C^\infty$ on $\Omega$. So the solution $RS + AF$ is of class $C^\infty$ on $\Omega$. It remains to show that $AF$ has a boundary value in currents sense on $\Omega$.\\
Let $(\varphi_j)_{j \in J}$ be a partition of the unit subordinated to a finite recovery $(U_j)_{j \in J}$ of $\bar{\Omega}$ by the open set of local coordinates.\\
We have $AF = \displaystyle {\sum_{j \in J} A \varphi_jF }$ with $A\varphi_j F$ have compact support in $U_ j$.
\begin{center}
$AF_{\vert \Omega} = \displaystyle {\sum_{j \in J} (A \varphi_j F)_{\vert \Omega}}$.
\end{center}
If $U_j \subset \Omega$, then $ A \varphi_j F$ is of class $C^\infty$ with compact support in $\Omega$, so $(A \varphi_j F)_{\vert \Omega}$ admits a boundary value in currents sense.\\
If $U_j \nsubseteq \Omega$ and $U_j \cap b \Omega \neq \emptyset$; let us show that $ A \varphi_j F$ admits a boundary value in currents sense.\\
Since $\varphi_j$ is with support in $U_j$ which is an open set of coordinate, so we are reduced to a bounded domain of $\mathbb{C}^n$. $ A \varphi_j F$ is of the same nature than the action of the newton kernel $E(x, y)$ on $\varphi_j F$. Since $F$ is order $m$ with compact support in $\bar{\Omega}$ which prolongs $f$, thus according to Theorem \ref{1}, $ A \varphi_j F$ admits a boundary value in currents sense.\\
Since $Card(J)< \infty$, so $\displaystyle {\sum_{j \in J} (A \varphi_j F)_{\vert \Omega}}$ admits a boundary value in currents sense. Thus we have
\begin{center}
$\tilde{H}^{n,r}(\Omega) = 0$.
\end{center}
\end{preuve}
Taking into account Theorem \ref{1} and the results of the $\bar{\partial}$ for differential forms of class $C^\infty$ with  boundary value in currents sense in Corollary \ref{4}, the proof of Theorem \ref{zhou} can be made:
\begin{preuve} [theoreme \ref{zhou}] \item
Let $f$ be a $(p, q)$-form of class $C^\infty$, $d$-closed with boundary value in currents sense and defined on $\Omega$ for $1 \leq p + q \leq 2n - 1$. Then according to Theorem \ref{1}, there exists a $(p + q - 1)$-form $g$ of class $C^\infty$ with  boundary value in currents sense such that $dg = f$.\\
We do not lose in general considering that $g$ is broken down into a $(p - 1, q)$-form $g_1$ of class $C^\infty$ with boundary value in currents sense and in a $(p, q - 1)$-form $g_2$ of class $C^\infty$ with boundary value in currents sense.\\
We have $dg = d (g_1 + g_2) = dg_1 + dg_2 = f$.\\
Like $d = \partial + \bar{\partial}$, we have for reasons of bidegre $\partial g_2 = 0$ and $\bar{\partial} g_1 = 0$. According to Corollary \ref{4}, we have $g_1 = \bar{\partial}h_1$ and $g_2 = \partial h_2$ where $h_1$  and $h_2$ are a differential forms of class $C^\infty$ with boundary value in currents sense defined on $\Omega$.
We have $ f = \partial g_1 + \bar{\partial}g_2 = \partial \bar{\partial}h_1 + \bar{\partial} \partial h_2$\\
Now $\partial \bar{\partial} = - \bar{\partial}\partial$
\begin{center}
 $\partial \bar{\partial} h_1 - \partial \bar{\partial} h_2 = \partial \bar{\partial} (h_1-h_2)= f$.
 \end{center} 
Let $u = h_1 - h_2$, $u$ is a $(p - 1, q - 1)$-form of class $C^\infty$ with boundary value in currents sense defined on $\Omega$ such that $\partial \bar{\partial} u = f$.
\end{preuve}

Better, taking into account Theorem \ref{q}, we obtain the following result:
\begin{theorem} \label{sam}
Let $M$ be a complex analytic manifold of dimension $n$
and $\Omega \subset \subset M$ be a contractible domain with smooth boundary of class $C^\infty$. Suppose
that $\Omega$ is completely strictly $q$-convex for $0 \leq q \leq n-1$ and $H^j(b\Omega)$ is trivial for $1 \leq j \leq 2n - 2$. Then the equation  $\partial \bar{\partial} u = f$, where $f$ is a $(n, r)$-form of class $C^\infty$, $d$-closed with boundary value in currents sense for $n - q + 1 \leq r \leq n - 1$ admits an solution $u$ which is a $(n - 1, r - 1)$-form of class $C^\infty$ with boundary value in currents sense.
\end{theorem}
\begin{preuve}
Let $f$ be a $(n, r)$-form of class $C^\infty$, $d$-closed with boundary value in currents sense defined on $\Omega$. So According to the theorem \ref{1}, there exists a $(n + r - 1)$-form $g$ of class $C^\infty$ with boundary value in currents sense such that $dg = f$.\\
We do not lose in general considering that $g$ is broken down into a $(n - 1, r)$-form $g_1$ of class $C^\infty$ with boundary value in currents sense and in a $(n, r-1)$-form $g_2$ of class $C^\infty$ with boundary value in currents sense.\\
\begin{center}
$dg = d (g 1 + g 2) = dg 1 + dg 2 = f$.
\end{center}
Now $d = \partial + \bar{\partial}$, we have for reasons of bidegre $\partial g_2 = 0$ and $\bar{\partial}g_1 = 0$.\\
According to Theorem \ref{q}, we have $g_1 = \bar{\partial}h_1$ and $g_2 = \partial h_2$ where $h_1$ and $h_2$ are a differential forms class $C^\infty$ with boundary value in currents sense defined on $\Omega$.\\
We have $ f = \partial g_1 + \bar{\partial}g_2 = \partial \bar{\partial}h_1 + \bar{\partial} \partial h_2$\\
Now $\partial \bar{\partial} = - \bar{\partial}\partial$
\begin{center}
 $\partial \bar{\partial} h_1 - \partial \bar{\partial} h_2 = \partial \bar{\partial} (h_1-h_2)= f$.
 \end{center} 
Let $u = h_1 - h_2$, $u$ is a $(n - 1, r - 1)$-form of class $C^\infty$ with boundary value in currents sense defined on $\Omega$ such that $\partial \bar{\partial} u = f$. 
\end{preuve}
\section{ Solving $\partial \bar{\partial}$ for a differential forms with boundary value in currents sense in $M \setminus \bar{D}$  where $D$ is a contractible completely strictly pseudoconvex domain}
Let $\Omega$ be a domain of a differentiable manifold $M$ of dimension $n$. In this part, it's about giving the analogous of theorem \ref{zhou} for $\Omega = M \setminus \bar{D}$, where $D \subset \subset M$ is a contractible completely strictly pseudoconvex domain and verifying $H^j(bD)$ trivial for $1 \leq j \leq n - 2$. For this, started by giving the following results:
\begin{theorem}\label{aj}
Let $X$ be a complex analytic manifold of dimension $n$ and $D$ be a domain on board $C^\infty$ strictly $q$-concave. For all $\xi \in bD$, there is a neighborhood $\theta$ of $\xi$, such that for any domain $D_1$ with on board $C^\infty$ sufficiently near to $D$ in the sense of the topology $C^2$ and for all $f \in \tilde{D}_{D_1}^{0,r}(X) \cap \ker \bar{\partial}$, $\bar{\partial}$-exact in $D_1$, with $1 \leq r \leq q$, there exists a $(0, r-1)$-form $g$ of class $C^\infty$ in $D_1 \cap \theta$ with boundary value in currents sense on $bD_1 \cap \theta$ such as $\bar{\partial}g = f$ on $D_1 \cap \theta$.
\end{theorem}
\begin{preuve}
It is identical to that of lemma $4.3$ of $[9]$, we build $\Delta$, $\Delta^{'}$ and $\Omega$ as in the proof of lemma $4.3$ of $[9]$ and replace the extensible current $\check{T}$ by the current $[f]$ and the theorem $2$ of $[9]$ by the theorem $4.1$ of $[7]$.
\end{preuve}
As a consequence of this theorem, we have the following corollary:
\begin{corollary} \label{Gn}
Let $X$ be a Stein manifold of dimension $n$. Let $\Omega \subset X$ such as $X$ is an extension $(n-1)$-concave of $\Omega$. Let $f$ be a differential form of bidegre $(0, r)$ on $\Omega$ and $\bar{\partial}$-closed with boundary value in currents sense for $1 \leq r \leq n - 2$. There is a $(0, r - 1)$-form $g$ defined on $\Omega$ with boundary value in currents sense such as
$\bar{\partial} g = f$.
\end{corollary}
\begin{preuve}
It's about showing that $\tilde{H}^{0,r}(\Omega) = 0$ for $1 \leq r \leq n-2$.\\
Since $\Omega$ is concave, we have $\check{O}_\Omega = O_\Omega$. So $H^r(\Omega,\check{O}_\Omega)\simeq H^r(\Omega,O_\Omega)$.\\
Consider the following sequence:
\begin{equation} \label{we}
0 \rightarrow \check{O}_\Omega \rightarrow \mathcal{F}^{0,0}(\Omega) \rightarrow \mathcal{F}^{0,1}(\Omega) \rightarrow  \cdots \rightarrow \mathcal{F}^{0,n-1}(\Omega) \rightarrow \bar{\partial}\mathcal{F}^{0,n-1}(\Omega) \rightarrow 0.
\end{equation}
According to Theorem \ref{aj}, we know how to solve locally the $\bar{\partial}$ for a differential forms
with boundary value in currents sense on $\Omega$. So the sequence \ref{we} is exact, and since the sheaf $\mathcal{F}^{0,r}(\Omega)$ is fine as a sheaf of modules on a sheaf of rings of class $C^\infty$, then the differential complex $(\mathcal{F}^{0,\bullet}(\Omega), \bar{\partial})$ of the differential forms defined on $\Omega$ with boundary value in currents sense is an acyclic resolution of the sheaf $\check{O}_\Omega$. By
hence for $0 \leq r \leq n - 2$, we have the following functorial isomorphism:\\
\begin{center}
 $H^r(\Omega,\check{O}_\Omega)\simeq \tilde{H}^{0,r}(\Omega):= \frac{ker(\bar{\partial}:\mathbf{E}^{0,r}(\Omega) \rightarrow \mathbf{E}^{0,r+1}(\Omega))}{Im(\bar{\partial}:\mathbf{E}^{0,r-1}(\Omega) \rightarrow \mathbf{E}^{0,r}(\Omega))}$
 \end{center}
where $\mathbf{E}^{0,r}(\Omega):= \Gamma(\Omega,\mathcal{F}^{0,r}(\Omega))$ are a sections on $\Omega$ of the sheaf $\mathcal{F}^{0,r}$.\\
So $H^r(\Omega,O_\Omega)\simeq \tilde{H}^{0,r}(\Omega)$. According to the isomorphism of Dolbeault, we have $H^r(\Omega,O_\Omega)= H^{0,r}(\Omega)$. Since $X$ is an extension $(n - 1)$-concave of
$\Omega$, we have by invariance of the cohomology: $H^{0,r}(\Omega) = H^{0,r}(X) = 0$ for $1 \leq r \leq n-1$.\\ 
So
\begin{center}
  $\tilde{H}^{0,r}(\Omega) = 0$, for $1 \leq r \leq n-2$.
  \end{center}
\end{preuve}
By taking inspiration from the proof of corollary \ref{Gn}, we obtain in the following theorem the global analogous of an result obtained in $[7]$.
\begin{theorem} \label{djilo}
Let $M$ be a complex analytic manifold of dimension $n$ and $D \subset \subset M$ be a completely strictly pseudoconvex domain with smooth boundary such as $M$ is an extension $(n-1)$-convex of $D$. Let's ask $\Omega = M \setminus \bar{D}$. Let $f$ be a differential form of bidegre $(0, r)$ defined on $\Omega$ of class $C^\infty$, $\bar{\partial}$-closed with boundary value in currents sense for $1 \leq r \leq n - 2$. There exists a $(0, r-1)$-form $g$ defined on $\Omega$ of class $C^\infty$ with boundary value in currents sense such as $\bar{\partial}g = f$.
\end{theorem}
\begin{preuve}
It's about showing that $\tilde{H}^{0,r}(\Omega) = 0$ for $1 \leq r \leq n-2$.\\
Since $\Omega$ is concave, we have $\check{O}_\Omega = O_\Omega$.  So $H^r(\Omega,\check{O}_\Omega)\simeq H^r(\Omega,O_\Omega)$.\\
Consider the following sequence:
\begin{equation} \label{wa}
0 \rightarrow \check{O}_\Omega \rightarrow \mathcal{F}^{0,0}(\Omega) \rightarrow \mathcal{F}^{0,1}(\Omega) \rightarrow  \cdots \rightarrow \mathcal{F}^{0,n-2}(\Omega) \rightarrow \bar{\partial}\mathcal{F}^{0,n-2}(\Omega) \rightarrow 0.
\end{equation}
According to Theorem \ref{aj}, we know how to solve locally the $\bar{\partial}$ for differential forms
with boundary value in currents sense on $\Omega$. So the sequence \ref{wa} is exact, and since the sheaf $\mathcal{F}^{0,r}(\Omega)$ is fine as a sheaf of modules on a sheaf of rings of class $C^\infty$, then the differential complex $(\mathcal{F}^{0,\bullet}(\Omega), \bar{\partial})$ of the differential forms defined on $\Omega$ with boundary value in currents sense is an acyclic resolution of the sheaf $\check{O}_\Omega$. By hence for $0 \leq r \leq n - 2$, we have the following functorial isomorphism:
\begin{center}
$H^r(\Omega,\check{O}_\Omega)\simeq \tilde{H}^{0,r}(\Omega)$
\end{center}

So $H^r(\Omega,O_\Omega)\simeq \tilde{H}^{0,r}(\Omega)$. According to the isomorphism of Dolbeault, we have $H^r(\Omega,O_\Omega)= H^{0,r}(\Omega)$. Since $M$ is an extension $(n-1)$-concave of $\Omega$, we have by invariance of cohomology: $H^{0,r}(\Omega) = H^{0,r}(M)$ for $0 \leq r \leq n-2$.\\
According to $[3]$ $H^{0,r}(D) \simeq H^{0,r}(M) = 0$ for $1 \leq r \leq n$.\\
So\\
\begin{center}
  $\tilde{H}^{0,r}(\Omega) = 0$, pour $1 \leq r \leq n-2$.
  \end{center}
\end{preuve}
More generally, if $0 < q < n - 1$, we have the following theorem:
\begin{theorem} \label{djilo}
Let $M$ be a complex analytic manifold of dimension $n$ and $D \subset \subset M$ a completely strictly pseudoconvex domain with smooth boundary such as $M$ is an extension $q$-convex of $D$ for $q \geq \frac{n+1}{2}$. Let's ask $\Omega = M \setminus \bar{D}$. Let $f$ be a differential form of bidegre $(0, r)$ defined on $\Omega$ of class $C^\infty$, $\bar{\partial}$-closed with boundary value in currents sense for $1 \leq n - q \leq r \leq q - 1$. There exists a $(0, r - 1)$-differential form $g$ defined on $\Omega$ of class $C^\infty$ with boundary value in currents sense such that $\bar{\partial}g = f$. 
\end{theorem}
\begin{preuve}
It's about showing that $\tilde{H}^{0,r}(\Omega) = 0$ for $n-q \leq r \leq q-1$\\
Since $\Omega$ is concave, we have $\check{O}_\Omega = O_\Omega$. So $H^r(\Omega,\check{O}_\Omega)\simeq H^r(\Omega,O_\Omega)$.\\
Consider the following sequence: 
\begin{equation} \label{wi}
0 \rightarrow \check{O}_\Omega \rightarrow \mathcal{F}^{0,0}(\Omega) \rightarrow \mathcal{F}^{0,1}(\Omega) \rightarrow  \cdots \rightarrow \mathcal{F}^{0,n-2}(\Omega) \rightarrow \bar{\partial}\mathcal{F}^{0,n-2}(\Omega) \rightarrow 0
\end{equation}
According to Theorem \ref{aj}, we know how to solve locally the problem of $\bar{\partial}$ for differential forms with boundary value in currents sense on $\Omega$. So the sequence \ref{wi} is exact, and since the sheaf $\mathcal{F}^{0,r}(\Omega)$ is fine as a sheaf of modules on a sheaf of rings of class $C^\infty$, then the differential complex $(\mathcal{F}^{0,\bullet}(\Omega), \bar{\partial})$ of the differential forms  defined on $\Omega$ with boundary value in currents sense is an acyclic resolution of the sheaf $\check{O}_\Omega$. By hence for $0 \leq r \leq n - 2$, we have the following functorial isomorphism:
\begin{center}
$H^r(\Omega,\check{O}_\Omega)\simeq \tilde{H}^{0,r}(\Omega)$
\end{center}

So $H^r(\Omega,O_\Omega)\simeq \tilde{H}^{0,r}(\Omega)$. According to the isomorphism of Dolbeault, we have $H^r(\Omega,O_\Omega)= H^{0,r}(\Omega)$. Since $M$ is an extension $q$-concave of $\Omega$, we have
by cohomology invariance: $H^{0,r}(\Omega) = H^{0,r}(M)$ for $0 \leq r \leq q-1$.\\ 
From $[3]$ $H^{0,r}(D) \simeq H^{0,r}(M) = 0$ for $n-q \leq r \leq n$.\\
So
\begin{center}
  $\tilde{H}^{0,r}(\Omega) = 0$, pour $n-q \leq r \leq q-1$.
  \end{center}
\end{preuve}
To prove theorem \ref{aff}, we need the following lemma:
\begin{lemma} \label{aff1}
Let $M$ be a differentiable manifold of dimension $n$, $D \subset \subset M$ be a ontractible domain with smooth boundary of class $C^\infty$ such that $M$ is an contractible extension of $D$ with $H^j(bD)$  trivial  for $1 \leq j \leq n-2$. Let $\Omega = M \setminus \bar{D}$. So $\stackrel{\circ}{\bar{\Omega}} = \Omega$ and $f$ is a $r$-form of class $C^\infty$, $d$-closed with boundary value in currents sense on $\Omega$ for $0 \leq r \leq n$, then there exists a $(r-1)$-form $g$ of class $C^\infty$ defined on
$\Omega$ with boundary value in currents sense such that $dg = f$.
\end{lemma}
\begin{preuve}
Let $f$ be a differential form of class $C^\infty$, $d$-closed with boundary value in currents sense on $\Omega$, according to Lemma $4.1$ of $[7]$, $[f]$ is a extensible current of finite order. According to $[2]$,  $\check{H}^r(\Omega) = 0$ for $1 \leq r \leq n-1$. 
So there exists a $(r-1)$-extensible current $u$ defined on $\Omega$ such that $du = f$.\\
Let $S$ be an extension of $u$ on $M$. Let $dS = F$, $F$ is an extension
of the current $[f]$ and therefore of finite order $m$. We have
\begin{center}
$S=RS +AdS + dAS$.
\end{center}
Let $\check{S} = RS +AdS$, $RS$ is an differential form of class $C^\infty$ on $M$, therefore $RS_{\vert \Omega}$ admits a boundary value in currents sense. Since $A$ does not increase the singular support so $AdS_{\vert \Omega}$ is of class $C^\infty$ and admits as in the proof of Theorem \ref{1} a boundary value in currents sense.\\
Thus $\check{S}_{\vert \Omega}$ admits a boundary value in currents sense and $d \check{S}_{\vert \Omega} = f$.
\end{preuve}
\begin{preuve} [Theoreme \ref{aff}] \item 
Let $f$ be a $(p, q)$-form of class $C^\infty$, $d$-closed with boundary value in currents sense on $\Omega$. According to Lemma \ref{aff1}, there exists a $(p + q - 1)$-form $h$ of class $C^\infty$ with boundary value in currents sense such that $dh = f$. We have $h = h_1 + h_2$ where $h_1$ and $h_2$ are respectively a $(P - 1, q)$-form class $C^\infty$ with boundary value in currents sense and a $(p, q - 1)$-form of class $C^\infty$ with boundary value in currents sense. We have
\begin{center}
$dh= dh_1 +dh_2=f$.
\end{center}
Since $d= \partial + \bar{\partial}$ and for reasons of bidegre we have $\partial h_2 =0$ and $\bar{\partial} h_1 =0$. According to Theorem \ref{djilo}, $h_1 =\bar{\partial}g_1$ and $h_2 = \partial g_2$ where $g_1$ and $g_2$ are a differentials forms of class $C^\infty$ with boundary value in currents sense defined on $\Omega$. We have
\begin{center}
$f = \partial h_1+ \bar{\partial}h_2= \partial \bar{\partial}g_1 +\bar{\partial}\partial g_2 = \partial \bar{\partial}(g_1 - g_2)$.
\end{center}
Let $u= g_1 - g_2$, $u$ is a $(p - 1, q - 1)$-form of class $C^\infty$ with boundary value in currents sense defined on $\Omega$ such that $\partial \bar{\partial} u = f$.
\end{preuve}
\section*{References}
 

\begin{enumerate}
 \bibitem [1]{1} E. BODIAN, D. DIALLO, S. SAMBOU: R{\'e}solution du $\partial \bar{\partial}$ pour les courants prolongeable d{\'e}finis sur un domaine compl{\'e}tement strictement pseudoconvexe d'une vari{\'e}t{\'e} analytique complexe, Imhotep Mathematical Journal, Volume $3$, Num{\'e}ro $1$, $(2018)$, pp. $1-4$.
\vspace{0,3cm}
\bibitem [2]{2} E. BODIAN, I. HAMIDINE, S. SAMBOU: R{\'e}solution du $\partial\bar\partial$ pour les courants prolongeables d{\'e}finis dans un anneaux, \url{arXiv:2250592v1}.
\vspace{0.3cm}
\bibitem [3]{3} G. M. HENKIN, J. LEITERER: Andreotti-Grauert theory by integrals formulas, Birkh{\"a}user, $1986$.
\vspace{0.3cm}
\bibitem [4]{4} C. LAURENT-THIEBAUT: Th{\'e}orie des fonctions holomorphes de plusieurs variables, Inter-{\'E}ditions et CNRS {\'E}ditions, $1997$.
\vspace{0.3cm} 
\bibitem [5]{5} S. LOJACIEWIECZ, G. TOMASSINI: Valeurs au bord des formes holomorphes, in Several Complex Variables (P. Scuola. Norm. Sup. Pisa,éd), Cortona, $197677$, $1978$, p. $222-246$.
\vspace{0.3cm}
\bibitem [6]{6} A. MARTINEAU: Distribution et valeurs au bord des fonctions holomorphes, Strasbourg RCP $25$, $1966$.  
\vspace{0.3cm}
 \bibitem [7]{7} S. SAMBOU, S. SAMBOU: R{\'e}solution du $\partial \bar{\partial}$ pour les formes diff{\'e}rentielles ayant une valeur au bord au sens des courants dans un domaine strictement pseudoconvexe, Annales Math{\'e}matiques Blaise Pascal, vol.$25$, $n^\circ$ $2$ ($2018$), p. $315-326$. 
\vspace{0.3cm}
\bibitem [8]{8} S. SAMBOU: R{\'e}solution du $\bar\partial$ pour les courants prolongeables, Math. Nachrichten, $235$ , $179-190$, $2002$.
\vspace{0,3cm}
\bibitem [9]{9} S. SAMBOU: R{\'e}solution du $\bar\partial$ pour les courants prolongeables d{\'e}finis dans un anneau, Annales de la Facult{\'e} des sciences de Toulouse: Math{\'e}matiques, S{\'e}r. $6$, $11$ no. $1$, $2002$, $105-129$.
\end{enumerate}
\end{document}